# Diversity Dependent Uncertainty Management for Hydrocarbon Stimulation in Uncertain Heterogeneous Reservoir with Improved Efficiency

C. Cheng[*]

*Department of Chemical and Petroleum Engineering, University of Pittsburgh, Pittsburgh, PA, USA*

[*]*Corresponding Author:* chc203@pitt.edu

## Abstract

The history of oil and gas well stimulation through hydraulic fracturing is characterized by a pursuit of optimal designs tailored to reservoir properties. However, as with many engineering systems, the impact of variability and uncertainty (in this case of reservoir properties) is not accounted in these deterministic approaches. We propose an expansion of the principle of design diversity wherein the focus becomes developing optimal portfolio combinations of multiple designs rather than repeated application of a single design. This approach can substantially increase efficiency and decrease risk when designing systems where outcome desirability is definable as a smoothly-varying efficiency and where risk is related to the variance of outcomes tied to uncertain spatio-temporal variability of the design environment. For the case of oil/gas wells, we demonstrate diverse portfolios of designs can improve stimulation efficiency by nearly a factor of four and cut risk associated with variability and uncertainty of outcomes by over 80%. Hydraulic fracturing optimization is therefore fundamentally dependent upon uncertainty, and a change of approach toward the pursuit of design diversity can have a profound impact on efficiency and risk, and in turn cost and environmental impacts, in oil and gas development.

# 1. Introduction

Hydraulic fracturing (HF), a well stimulation technique in which rock is fractured by pressurized liquid, plays and ever-increasing role in energy production from oil and gas (hydrocarbon) resources. As of 2019, hydraulic fracturing enables about 75% U.S. dry natural gas production (1) and about 63% of total U.S. crude oil production (2). However, with this prominence comes a growing need to reduce negative environmental impacts (3,4). One metric is to consider environmental footprints caused by producing a unit of hydrocarbon, giving a type of efficiency of use of resources and impacts associated with hydrocarbon production (5,6). However, it is a major challenge to optimize this efficiency because it is, to a large degree, controlled by the variability and uncertainty of reservoir conditions, resulting in substantial uncertainty in HF stimulation and eventual production. The variability is naturally-possessed due to the sedimentation, tectonic movement, and other factors impacting formation of reservoir rocks (7,8). The variability is also human-induced due to variation in both design and execution of the processes involved in drilling, completing, and producing a well (9-12). Besides being variable, limitations of existing reservoir characterization technologies ensure that variability partially or fully transfers to become uncertainty because it is still impossible to measure formation variation accurately and economically for every well and throughout the reservoir away from the well.

Well stimulation technology trends are often driven by mitigating this variability and uncertainty because they represent risk associated with developing a given resource. Over recent years there is growing use of a method known as extreme limited entry (EXL), where the principle is to dominate the global pressure by the pressure drop through the well perforations by creating fewer and/or smaller perforation holes in each perforation cluster along the well (each perforation cluster represents the fluid entry point to an intended hydraulic fracture). Because perforations are man-made and, at least to a major degree, controllable, EXL is effective because it allows an engineered factor to determine distribution of fluid to each fracture within a stimulation stage (typically 3-6 fractures that are intended to be created simultaneously – see (13) for a description of multistage hydraulic fracturing of horizontal wells). Hence, variability in HF growth is less susceptible to be controlled by variable and uncertain distributions of stress and other rock properties along the well (14-16). However, there is a cost to this reduction in uncertainty, namely the need to significantly increase the pumping power required to sustain the increased pumping pressure while maintaining sufficient injection rate to grow the fractures within an economically and operationally viable timeframe. Hence, the economics of the well are negatively impacted and the environmental footprint is increased by the need for more pumps with larger fuel consumption.

Extreme limited entry (EXL) is one example illustrating that a drive to reduce variability and uncertainty of outcomes increases monetary and environmental cost associated with producing each unit of resource. But it is just one indicator of a broad trend towards larger, more economically and environmentally costly treatments that may be neither increasingly efficient nor increasingly effective in reducing variability of outcomes (i.e. risk). What is clear is that this trend can increase pollution (17-19), waste resources (20,21), and create growing risk and cost factors associated with water management and the potential for stimulation and/or wastewater injection to induce problematic seismicity (22-24).

Besides its importance as a technology, oil/gas well stimulation represents a useful archetype of a particular class of engineering design problem. Here we identify the essential nature defining this class of problem and present a new direction for increasing efficiency and reducing risk that is inspired by Modern Portfolio Theory (MPT) and Capital Market Link (CML), which are long-established methods in economics and finance for optimally balancing risk versus return in markets (25,26). The concept is that, a diversified portfolio of designs can reduce risk relative to repeated application of the same design in multiple wells. By using our validated rapidly-computing hydraulic fracture simulator Cheng and Bunger (2019), it is possible to both explore many different possible stimulation designs by varying factors like fracture spacing, injection rate and time, number of fractures per stage, number of perforation holes per cluster, and so forth. Furthermore, it is possible to generate multiple realizations of each design with rock properties and in-situ stresses drawn from stochastic distributions that represent the variability of these variable and uncertain inputs. The result is that each candidate design can be associated with an expected value of the return, which is in this case illustrated as an efficiency, defined as created fracture surface area divided by the total expended pumping energy. Each candidate design is also associated with risk, here illustrated by a combination of the variability in the total created fracture surface and the uniformity with which the stimulation penetrates the reservoir rock.

The framework inspired by MPT leads to an efficient frontier that emerges in this return versus risk space. The risk-free asset appears in CML inspire the role of low efficiency, low risk methods like EXL. EXL plays a role analogous to a government bond – it is typically lower return but lower risk. Such a "risk-free rate" can be combined with higher efficiency methods in an optimal way creating a return versus risk relationship tailored to a given situation. In this manner, a striking conclusion is shown to emerge: Rather than pursuing single optimal design, a novel approach of creating an optimized portfolio combination of high efficiency portfolio and low risk portfolio is invoked and identified as a promising way to manage shown capable of reducing risk while avoiding the loss in efficiency associated with high energy input methods such as extreme limited entry. The underlying concept is that the designs in high efficiency portfolio are tailored to higher returns while the designs in low risk portfolio are preferable to reduce risk, similar to the nature of various investments in a financial portfolio. And, just as in an optimally-diversified investment portfolio, mixing the proportion of wells stimulated with each of several chosen designs allows leveraging the advantages of each design rather than being beholden to the characteristics of a single design applied repetitiously across all wells.

The organization of the paper begins with a background description of Modern Portfolio Theory and Capital Market Link, describing the adaptations required for it to be applied for designing resource recovery with uncertain and variable reservoir conditions. Next, the methods are laid out with particular focus on the definitions used for return and risk. For the specific case of hydraulic fracturing stimulation of oil and gas wells, quantifying return and risk to identify portfolios of designs and explore the potential benefit of diverse designs is requires computationally-intensive analysis that is enabled by a very rapidly computing reduced order model of hydraulic fracture growth (27). The paper then describes the specifics of method deployment in order to firstly identify several individual strongly-performing HF designs and then to seek optimal combinations of these designs that creating high efficiency portfolios and low risk portfolios for portfolio combination, which optimally balance risk and return. Finally, we address the essential features of a class of design under uncertainty, for which the oil/gas well stimulation problem provides a archetype, and thereby providing a pathway for the use of diverse portfolios of design to increase efficiency and reduce risk across application areas ranging from renewable energy to agriculture to epidemioly.

## 2. Methods

## 2.1 Background

Encountering formation variability during hydraulic fracture (HF) stimulation of hydrocarbon wells is inevitable. Most efforts to optimize HF stimulation take a deterministic approach that does not account for variability. And, more subtly but in many ways more importantly, the existing optimization under uncertainty and variability takes a single well approach when, in fact, variability also exists between wells. Often there are hundreds to many thousands of wells drilled into the same reservoir, and an optimized design for a single well will perform differently from well to well. Maximizing return (i.e. efficiency of resource use and recovery) while minimizing risk (i.e. variability of outcomes and non-uniformity of reservoir stimulation) across many wells in a variable reservoir remains poorly understood and in need of new frameworks and paradigms bringing useful formalism to the optimization of well designs. While the problem of HF stimulation has some unique features, optimizing the risk/return tradeoff in variable and uncertain contexts is commonly encountered in other fields, especially economics and finance. Here we adopt two seminal works from these fields to provide a new framework for optimizing the risk/return tradeoff in HF stimulation.

The first seminal work from the field of finance is Modern Portfolio Theory, first proposed in 1952 as a quantitative framework to maximize the return of investments while minimizing risk by diversifying assets (25). The basic approach begins with a step wherein various candidate investment portfolios are assembled from assets that are mixed with various allocations. Each candidate portfolio is associated with a particular expected value of the return on investment as well as a given level of risk, i.e., the variance on the rate of return. Generating a cross plot of return versus risk for each of the candidate portfolios leads to a cloud of points, illustrated in Figure 1. In this way, the multiobjective optimization problem of maximizing return while minimizing risk is formally recognized. Using terminology of multiobjective optimization (28), this cloud of points represents the so-called feasibility set for a given collection of assets representing the return and risk for each mixture of the assets. The upper left boundary of the feasibility set defines a subset of portfolios (the Pareto set) providing the minimal risk for a given return. For all cases belonging to this Pareto set, the expected return is maximized within a given level of risk.

While MPT provides an elegant way to identify optimal performance, the specific risk/return ratios that can be accessed are relatively limited; only those ratios which happen to fall along the Pareto frontier can be obtained for a given mix of assets. If the level of risk is still unacceptably high, then the so-called Capital Market Link (CML) provides an approach to mix higher risk higher return assets (traditionally stocks) with lower risk lower return assets of a different type (traditionally government bonds or savings accounts) in order to achieve an acceptable risk level (26). The method entails plotting the "Risk Free Rate" (i.e. insured savings rate) in the return/risk cross plot (Figure 1), and then drawing a line from this rate that is tangent to the Pareto frontier. Hence, the tangency point represents a portfolio with 100% of holdings in a diversified choice of investments drawn from the Pareto set. At the other end, the Risk Free Rate point represents investing solely in insured savings. The points on the CML line between the tangency point and the Risk Free Rate represent various proportions of investments allocated to the tangency portfolio and to insured savings. The principle of CML is that the expected return is maximized for a given level of risk for all choices along this line.

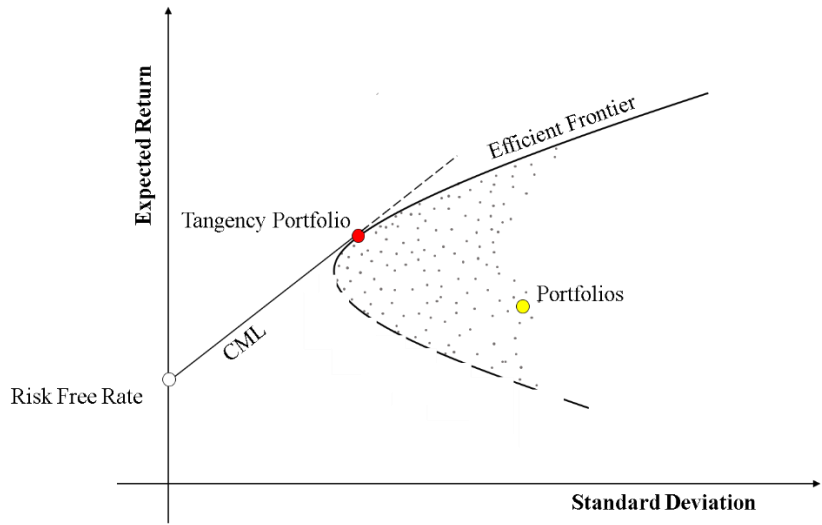

**Fig. 1.** Illustration of MPT and CML. Each choice of diversified mixture of assets in a portfolio results in a given expected value of the return and standard deviation of the return (i.e. risk), illustrated for a single case by the yellow dot. The black points are a feasibility set of portfolios with various mixtures of assets. The solid curve is the Pareto efficient frontier. The Risk Free Rate indicates the insured rate of return (i.e. savings rate). The straight line connecting the Risk Free Rate to the Tangency Portfolio (red dot) is referred to as the CML.

Developing MPT introduced the formal concept of diversification into finance and investment and, as such, received the 1990 Nobel Prize in Economics (25). However, if each possible HF design corresponds to a unique balance of risk and return in an uncertain and variable environment, then we hypothesize that diversification of design into a portfolio should provide preferable risk/reward tradeoff and MPT should provide a useful framework for evaluating this tradeoff. In other words, a portfolio constituted by distinctive designs with corresponding return and risk provides a chance to mitigate negative impacts of variability and uncertainty. Furthermore, if Extreme Limited Entry (EXL) provides an approach to low return low efficiency but minimal risk, then CML gives a rational way to mix EXL with other designs to achieve a desirable balance of efficiency and risk.

Clearly there will be some differences between investing and well stimulation. Most notably, the prices of assets are generally determined by the action of the same or related markets, with the result being readily-quantified covariance among assets that perform similarly in those markets. In well stimulation, the designs act upon different wells which, themselves, are variable and uncertain. As a result of such differences, definitions of return and risk will need to be developed that are definable within an HF simulator.

In the following, energy efficiency is introduced as a metric of the return of HF designs. Efficiency is calculated as the energy accounting for fracturing as a proportion of total input energy. The premise is that one of the ultimate and mechanically-quantifiable outcome of HF stimulation is generation of fracture surface area. This surface area, if appropriately propped open by injection of granular "proppant" (29), is often taken to be directly related to potential of a stimulation to generate increased hydrocarbon production (29). The total energy, on the other hand, determines the requirement of fluid volume and pumping power, that is, resources that must be expended in order to carry out HF stimulation. So, for a given HF design, higher energy efficiency means creating more fracture area per given unit of inputs, in other words, higher resource usage efficiency.

A given HF design not only corresponds to an expected value of the efficiency (return), but also to a variation of energy efficiency from one fracture to the next, which is defined as energy variability. Lower energy variability promotes more evenly distributed stimulation within the reservoir, meaning more evenly propped fractured area and reduced tendency to leave resource in place. And, besides spatial variance of stimulation along a given well, variance can also occur owing to the uncertainty and variability of rock properties encountered from one well to another. Reducing this variance promotes outcomes of stimulations across a set of many wells that can be more precisely anticipated. Therefore, spatial variability within each well and efficiency variability from well to well comprise variation of outcomes of HF stimulation that are driven by variable and uncertain inputs, in particular, reservoir rock properties and stress conditions. In the following, we develop a measure of risk that combines both of these forms of variability.

Once defined for well stimulation, the efficiency and risk are then predicted by carrying out simulations of many designs with varied formation parameters using a rapidly-computing model (27). This reduced-order model has been previously validated with several high-fidelity models (16,30). Model results are then used to demonstrate the risk and return of each design, including portfolios comprised of several designs applied in various proportion, like a diversified portfolio of assets in finance. In so doing, the adaptation a new framework modified on o MPT and CML to engineering design under uncertainty also introduces an expansion of the concept of design diversity, originally introduced to reduce faults in computing architecture (cites), to encompass a new approach to design that eschews the tendency to look for a single "best" design and instead looks for efficient combinations of multiple designs that can increase return and decrease risk when applied to systems with uncertain characteristics.

## 2.2 Individual Candidate Selection

Applying an adaptation of MPT and CML to hydraulic fracturing requires several steps. First, in an analogy to picking the individual stocks that will be mixed into a portfolio, various designs are simulated to identify high efficiency (i.e. high fracture area per unit energy input) HF designs with distinctive input parameters, that is, with sufficient difference among the designs that combining them would result in a diversified portfolio. The model considers an array of $N$ fractures symmetrically distributed within one stage of length $Z$ as illustrated by Figure 2(c) and the spacing $h_k$, $k=1,..N-1$ between each of the fractures is such that

$$Z = \sum_{k=1}^{N-1} h_k \tag{1}$$

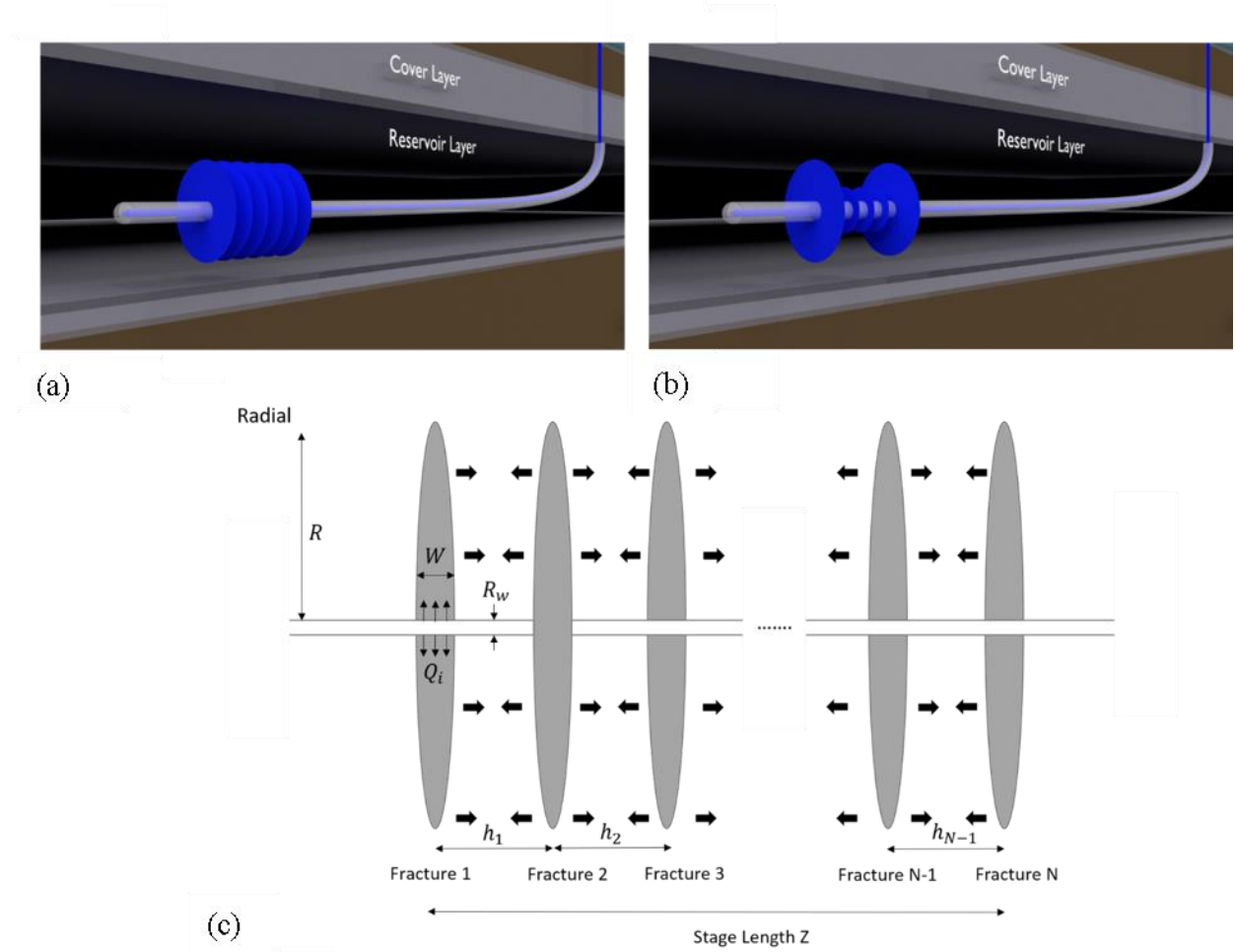

**Fig. 2.** Cross-sectional illustration of multiple, simultaneous HFs in one stage along a horizontal wellbore. (**a**) Ideal, uniform result, and (**b**) Result in which central fractures are suppressed. (**c**) Geometry of the multiple HF problem for $N$ HFs distributed within a stage of length $Z$ and with fracture spacing $h_k$. The arrows illustrate the interaction stresses between fractures. Figure adapted from (28).

Specifically, we consider two major classes of design, one with a uniform spacing between five fractures and the other with non-uniform spacings that allow five fractures to be non-uniformly (but symmetrically) placed in one stage. With many possibilities, for illustration the uniform spacing means spacing ratio $2h_1/Z=0.5$, where $h_1$ is the space between clusters and $Z$ is the stage length. For non-uniform spacing, $h_1$ is taken as the spacing between clusters 1 and 2 (and symmetrically between clusters 4 and 5) in a five-cluster array. The spacing ratio for non-uniform spacing varies from 0.25 to 0.5, where 0.5 corresponds to the limit of uniform spacing and 0.25 corresponds to a spacing between the outer sets of clusters that is three times smaller than the spacing between the inner clusters. The use of non-uniform spacing cases is inspired by past simulations (16,35) and laboratory experiments (36) showing that certain choices apparently balances the stresses exerted by fractures on one another – the so-called "stress shadow" effect – thereby mitigating the tendency of some fractures to be suppressed by the stresses generated by others Figure 2(b).

The search for designs considers 7200 candidates each for uniform and non-uniform spacing cases with varying values of: 1) Perforation pressure loss, i.e. more or fewer perforation holes per cluster leading perforation pressure loss between $10^1$ to $10^8$ Pa (details provided in the Supplement), 2) Injection rate $Q_o$ from $0.1 \text{ m}^3/\text{s}$ to $0.25 \text{ m}^3/\text{s}$, 3) Fluid viscosity $\mu$ from $10^{-3.5}$ to $10^1 \text{Pa.s}$, and 4) Stage length $Z$ from 20 to 120m. In reality there are more design parameters one could change including length of the productive horizontal portion of the well ("lateral length", here taken as 1000m), number of fractures per stage (here we restrict all designs to five fractures), and total pumping volume (here every well injects 14,400 $m^3$).

While rock properties will be varied later, at the initial point of selecting individual designs, realizations are generated for rock properties fixed as

$$C_{L0} = 3.24 \times 10^{-6} \text{ m/s}^{1/2}, E = 25 \text{ GPa}, \nu = 0.2 \tag{2}$$

$$\sigma_o = 34.47 \text{ MPa}, K_{IC} = 1 \text{ MPa} \cdot \text{m}^{1/2}$$

The leak-off coefficient $C_L$ is a function of both rock and fluid properties (37) and is therefore determined as $C_{L0}\sqrt{1\ Pa\ s/\mu}$ for a given fluid viscosity $m$, using $C_{L0}$ as a reference leak-off coefficient. Then, with these rock properties, each design is generated by random sampling of each design parameter utilizing Latin Hypercube Sampling to cover the parametric space efficiently. For each parameter selection, the inputs are provided to the rapidly-computing HF model C5Frac (see Supplement for details). The energy efficiency of each design is calculated through the ratio of the energy associated with rock breakage to the total energy input required by the system (see details in the Supplement, Section S2)

$$\varepsilon = \left[\frac{\pi(1-\nu^2)K_{IC}^2}{E}\sum_{i=1}^{N}R_i^2(T)\right] / \int_0^T p_{f(i)}(R_w, t)Q_i(t)dt \tag{3}$$

Here the model predicts, for each ($i^{th}$) fracture, the fluid pressure $p_{f(i)}$ and the inflow rate $Q_i$ to the $i^{th}$ fracture as well as the resulting fracture radius, $R_i$. Additionally $N$ is the selected number of fractures in one stage (fixed in this illustration as 5) and $T$ is the selected treating time (chosen based on the injection rate to ensure equal injected volume for each case). Although rock property variation will be introduced after this initial selection, variability of fracture growth still occurs within each stage due to stress interaction whereby some fractures are suppressed and other are favored ("stress shadow") Figure 2(b). The model-predicted energy variability $\mathcal{V}$ is calculated through the standard deviation of energy efficiency of each of the clusters, that is

$$\mathcal{V} = \sqrt{\frac{1}{N-1}\sum_{i=1}^{N}\left(\frac{\varepsilon_i - \bar{\varepsilon}}{\bar{\varepsilon}}\right)^2} \tag{4}$$

With the defined energy efficiency and energy variability, a Pareto frontier is obtained via a crossplot of energy variation versus energy efficiency for all cases (See Supplement for details, especially Figures S1a and S1b). Designs are then picked by uniformly gridding the Pareto frontier, and these chosen designs will eventually be combined in various ratios to build a portfolio seeking to leverage their high efficiency while using the portfolio combination to reduce variability. Additionally, for comparison, we pick six lower efficiency uniform spacing designs from below the 2.5*10-4 inflection point. These will eventually be shown as single designs, applied repeatedly rather than combined into portfolios.

## 2.3 Accounting for Rock Variability

While variation brought on deterministically due to stress shadow is already accounted for in step one, the second step accounts for rock properties that are both variable and uncertain. Here the uncertainty level is defined as the variation range of the formation parameters. Hence we are assuming there is no characterization effort that would reduce the uncertainty level to below the variability level – which is often a valid assumption for low-cost onshore oil and gas wells. Each simulation case therefore draws rock properties at random from uniform distributions limited by upper and lower values provided by uncertainty level. Clearly other choices of distribution are possible, but a full exploration of impact of distribution is not needed to illustrate the method. In this illustration, we apply same uncertainty level of 5% (+/- 2.5% from the values described in Equation 2) for in-situ tress, elasticity and toughness. Larger uncertainty of from $10^{-2.5\%}$ to $10^{2.5\%}$ of the base value in Equation 2 is taken for leak-off coefficient because the underlying rock permeability can vary over several orders of magnitude. Here the focus will initially be on a 5% variability illustration, but 10% and 20% will be discussed later and are also provided as additional supporting cases in the Supplement.

With the rock property distributions defined, the properties are then independently and randomly picked using Latin Hypercube Sampling. With these variable rock properties, the high efficiency designs, selected in step one, are simulated using C5Frac to compute the energy efficiency and energy variability. This process of selecting random rock properties is repeated 120 times for each design. The energy efficiency is then taken as the mean of the efficiencies computed for all 120 realizations. Furthermore, the resulting total risk for each case is taken as a sum of the standard deviation energy variability defined in Equation 4 and the standard deviation of the energy efficiency from these 120 simulations with varying rock properties, that is

$$\sigma = \sqrt{\frac{1}{N_R - 1} \sum_{n=1}^{N_R} \left(\frac{\mathcal{V}_n - \bar{\mathcal{V}}}{\bar{\mathcal{V}}}\right)^2} + \sqrt{\frac{1}{N_R - 1} \sum_{n=1}^{N_R} \left(\frac{\varepsilon_n - \bar{\varepsilon}}{\bar{\varepsilon}}\right)^2} \tag{5}$$

Here $N_R$ is the number of realizations of variable rock properties.

## 2.4 Assembling and Simulating Candidate Portfolios of Designs

The individual designs are next combined in various proportions in order to assemble design portfolios. Each portfolio is associated with its own efficiency and risk. The risk and return of portfolio are calculated though MPT, which provides a mathematical framework (25) to mix individual asset with associated specific risk $\sigma_i$, return $r_i$ and the weighting $\rho_i$ of component asset $i$ (that is, the proportion of asset "$i$" in the portfolio). For each portfolio, a vector is assembled of the expected values of the return, recalling this is the model-computed efficiency averaged over all sets of rock properties. This vector of expected return and a vector of the weighting ratios are given by

$$\boldsymbol{E}(r) = \begin{bmatrix} r_1 \\ \vdots \\ r_N \end{bmatrix}, \quad \boldsymbol{\rho} = \begin{bmatrix} \rho_1 \\ \vdots \\ \rho_N \end{bmatrix}, \quad \sum_{i=1}^{N} \rho_i = 1 \tag{6}$$

The expected return of portfolio is given by the inner product of these vectors, that is

$$E(r_p) = \boldsymbol{\rho}' \cdot \boldsymbol{E}(r) = [\rho_1 \cdots \rho_N] \begin{bmatrix} r_1 \\ \vdots \\ r_N \end{bmatrix} = \sum_{i=1}^{N} \rho_i r_i \qquad (7)$$

Again, following MPT, the variance of the portfolio is related to the covariance of its members. For investments, this captures the fact that stocks with highly correlated prices do not reduce a portfolio's risk through diversification as effectively as stocks with uncorrelated prices. Hence, the covariance for two random assets $i$ and $j$ in portfolio, each with standard deviation $\sigma$, is defined by

$$cov(i,j) = \begin{bmatrix} \sigma_{11} & \cdots & \sigma_{1N} \\ \vdots & \ddots & \vdots \\ \sigma_{N1} & \cdots & \sigma_{NN} \end{bmatrix}, \sigma_{ij} = \sigma_i \sigma_j c_{ij} \qquad (8)$$

where $c_{ij}$ gives the correlation coefficients between the returns on design $i$ and $j$, and noting that $c_{ij}$ is one when $i$ equals $j$. Recall here that the risk is computed based on HF simulations for varying rock properties via Equation 5. In MPT the covariance among different members of the portfolio is non-zero because they are all performing on the same market or under the same set of market drivers. However in HF, each well is acted upon by a single design. The rock properties for that well are randomly chosen in a range, reflecting that the naturally-possessed variability and *a priori* uncertainty precludes choosing a design from the portfolio that will be best suited to the well. Instead, that well receives one of the designs from the portfolio, chosen at random, and without correlation to its unknown properties. The result is a standard deviation of the portfolio that is a weighted sum of the standard deviation of each member of the portfolio, with the weights coming from the proportion of that design in the portfolio. The portfolio standard deviation is therefore equivalent to the MPT-derived standard deviation if the off-diagonal elements of the correlation coefficient matrix are taken as zero. The portfolio variance is thus given by a double inner product of the covariance with the ratio vector

$$\sigma_p^2 = \rho' cov(i,j) \rho = [\rho_1 \cdots \rho_N] \begin{bmatrix} \sigma_{11} & \cdots & \sigma_{1N} \\ \vdots & \ddots & \vdots \\ \sigma_{N1} & \cdots & \sigma_{NN} \end{bmatrix} \begin{bmatrix} \rho_1 \\ \vdots \\ \rho_N \end{bmatrix} = \sum_{i=1}^{N} \sum_{j=1}^{N} \rho_i \rho_j \sigma_{ij} \qquad (9)$$

The positive square root of portfolio variance is named as risk is henceforth denoted as $\mathcal{R}$. To approximate the feasibility set of portfolios with all combinations of mixing ratios *r*, we randomly choose 12,500 combinations using a Monte Carlo method.

## 2.5 Building Optimal Portfolio Combinations

More control of risk to match a given level of risk tolerance is enabled through an analogue approach to Capital Market Link (CML), which requires firstly defining a Risk Free Rate. However, because there is no perfect analog to a Risk Free Rate, the same mixture method is applied on Extreme Limited Entry (EXL) designs to build minimal risk portfolios, thereby enabling construction of multiple lines that are drawn tangent to the high efficiency portfolio. To do this, we set varied perforation pressure factors while fixing a constant fluid viscosity of $0.003 \text{Pa} \cdot \text{s}$, injection rate of $0.25 \text{ m}^3/\text{s}$, stage length of 20m, and uniform spacing of 5 m between each fracture. The level of perforation pressure factors leads to a uniformly-distributed range of perforation pressure drops from 5.5MPa to 100MPa. The upper end of this range is impractical due to exceeding of typical equipment pressure ratings, however the search range is extended this far in order to ensure the best approximation of the true minimum risk designs is devised using EXL, without regards for externalities like pressure ratings. Upon selecting these candidate individual designs, we use random formation variations on the low risk designs at same uncertain level as with high efficiency designs. Six designs with lowest risk are selected to constitute the minimum risk portfolio. The exact value of corresponding perforation factor of each design occupied in low risk portfolio are detailed in the Table 2.

Upon obtaining a low risk portfolio, multiple lines are projected from the low risk portfolios to the Pareto front of higher efficiency but higher risk portfolios. Finally, the optimal portfolio combinations are identified to lie on the newly predicted frontier, which is a line tangent to both the low risk portfolios frontier and the high efficiency portfolios frontier. The risk-return of a given portfolio combination is therefore lying along this tangent line, with its location determined by the proportion of the low risk portfolios contained within the portfolio combination. For example, a 50:50 mix of tangency low risk and tangency high efficiency portfolios will give a composite efficiency versus risk that is from the midpoint of the line drawn between these two portfolios.

## 3. Results

Applying the four steps of the method for the illustration case leads to combinations of selected designs into high efficiency and low risk portfolios. We begin by demonstrating these combinations and then showing how optimal pathways for additional risk reduction can be attained by combining high efficiency and minimal risk (EXL) designs. The details of creating feasibility sets of designs, as described in Sections S3.1-3.2, are included in the Supplementary materials, with the focus here on the behavior of the portfolio combinations of these designs and the observations that can be made of how these portfolio combinations impact the efficiency versus variability tradeoff.

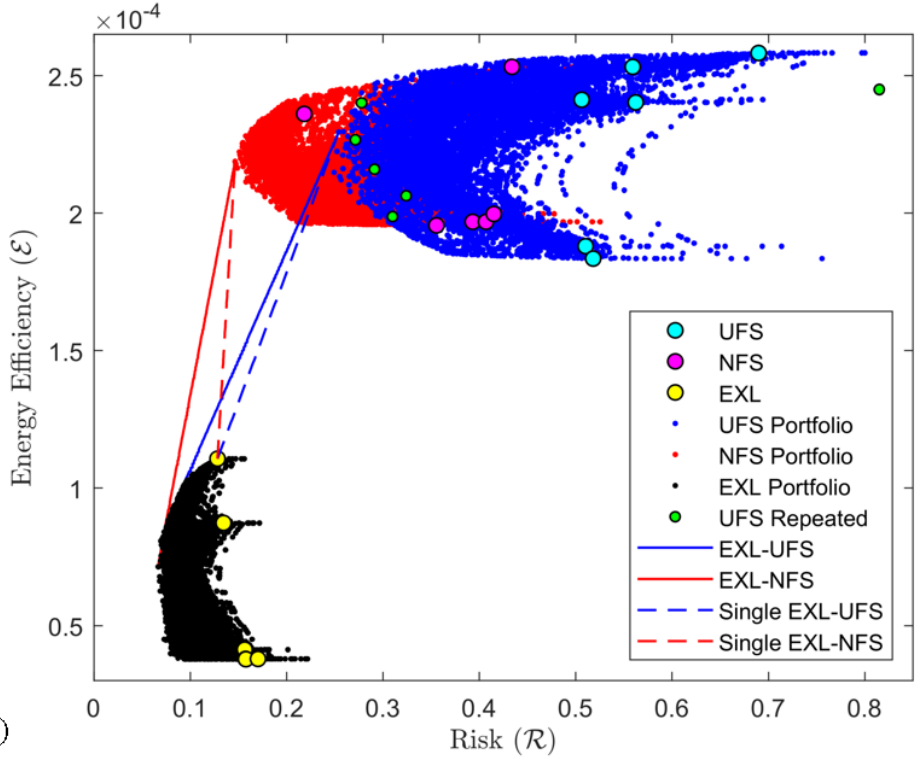

(a)

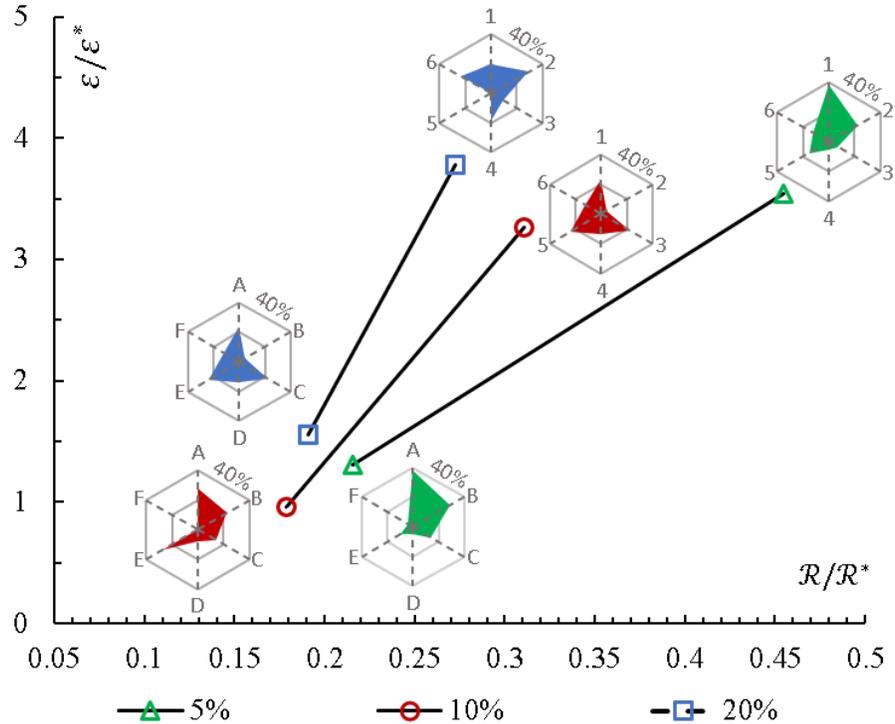

(b)
19

**Fig. 3.** Portfolios and optimal portfolio combinations. (a) The portfolios using uniform spacing (UFS) designs are plotted at right upper with highest risk as green color circle. The corresponding UFS portfolio is demonstrated by blue dots and the selected design for repetition is colored as cyan. The NFS portfolios are the red dots and repeated designs are colored as orange. The single extreme limited entry designs (EXL) are located in the lower left. Solid lines connecting the frontiers represent the CML for the UFS (red) and NFS (blue) cases. (b) Optimal portfolio combinations at different uncertainty levels in a plot of energy efficiency versus risk, scaled by a single design base case. Results are shown for 5%, 10%, and 20% uncertainty level in rock properties. The proportion of each design (detailed in Tables 1, and 2) are also given.

Through mixing candidate designs, feasibility sets of portfolios are obtained for uniform spacing (UFS), non-uniform spacing (NFS), and extreme limited entry (EXL) designs, as shown in Figure 3(a). These feasibility sets show that EXL designs, as expected, provide the lowest risk but also the lowest efficiency, owing to the high injection pressures required to sustain the large level of pressure loss through the perforations. It is therefore confirmed that EXL designs provide lowest efficiency and minimal risk. However, while this result is expected, what is more surprising is that the uniform and non-uniform spacing designs produce substantially different feasibility sets. Most strikingly, the variability of the outcomes, captured by the risk (Equation 9), is substantially different, with a large portion of the non-uniform spacing feasibility set showing lower risk compared to uniform spacing. Furthermore, the lower range of efficiencies for the uniform designs lie below the bottom of the feasibility set for the non-uniform designs, meaning that the "worst case" outcomes for non-uniform designs perform better than the worst case outcomes for uniform cases. It is worth noting, however, that this result is shown for 5% variability of the rock properties.

As the variability increases, in particular variability of the in situ stress, the random stresses dominate the stress interaction between the fractures and therefore advantages of non-uniform spacing can be eventually diminished (for this example, at around 20% variability of the in-situ stress there becomes little difference between uniform and non-uniform designs, see details in the Supplement).

While the potential risk reduction associated with non-uniform spacing is striking, perhaps the most important observation is the impact of using a portfolio of designs rather than repeating a single design on all wells. Examining the non-uniform example in Figure 3(a), the small red dots result from combinations of the single cases shown with large magenta dots. It is clear from comparing the feasibility set defined by the red dots to the magenta dots that a 40% reduction in risk is obtained from the portfolio Note this is obtained when comparing the portfolio to the repeated design with highest efficiency, that is, around $2.5*10^{-4}$. A similar level of risk reduction can be found upon examination of the single designs for uniform spacing in Figure 3(a) (large blue dots) to the feasibility set corresponding to combinations of the single designs (small green dots).

To reach even lower risk, six designs with lowest risk are chosen from extreme limited entry (EXL) and simulating random mixtures of these designs to assemble EXL portfolios, shown as the black points in Figure 3(a). Firstly we observe that, compared to the risk of a single EXL, the portfolios could contribute additional 15% risk reduction. Tangent lines (optimal lines showing results of mixing EXL and higher efficiency portfolios) are then drawn from the low risk EXL portfolios to the high efficiency portfolios. The line tangent to both the frontier of EXL portfolios and frontier of high efficiency portfolios gives an optimal line for which there is no portfolio locating above the line with higher efficiency at same risk, thereby confirming this line as an approximation of a Pareto frontier for mixtures of minimum risk and higher efficiency portfolios.

By this construction, the points presented on this line are optimal portfolio combinations representing the optimal mixtures of tangent high efficiency portfolio and corresponding tangency low risk portfolio for requested efficiency or limited risk. The change in efficiency is thus linearly related to the change in risk as the portfolio proportions between low risk and high efficiency portfolios vary. Because this is an approximation of a Pareto front, any other portfolios are predicted to induce lower efficiency at the same risk and/or higher risk at same efficiency. So, for example, for an energy efficiency of $2*10^{-4}$, an optimal portfolio combination could halve the risk of a non-uniform spacing (NFS) frontier portfolio by mixing NFS designs with EXL designs (i.e. running stimulating certain proportions or wells, or even stages within each well, with each design).

It is important to reiterate that the advantages of NFS compared to UFS designs are a consequence of the significant difference in the location of the Pareto frontier between NFS portfolio and UFS portfolio. However, as previously discussed, when the variability increases with higher uncertainty level, the interaction stress acted upon by optimized spacing will be overwhelmed by the in-situ stress variability. The spacing effect will be relatively lessened until the frontier of NFS portfolio frontier become similar to UFS portfolio. Thus, the nature of the optimal high efficiency portfolio and the optimal line connecting this portfolio to the minimum risk EXL portfolio turns out to be substantially impacted by the level of variability of rock properties. This clearly implies that, as one might expect intuitively, highly variable reservoirs will be optimally developed by different portfolios compared to less variable reservoirs.

To see more clearly the impact of variability, three levels are chosen: 5%, 10% and 20% (the corresponding variation is detailed in Supplement Section S3.3). Three optimal portfolio combinations are developed corresponding to each of these three uncertainty levels. It is most illuminating to examine the relative improvement (i.e. increase in efficiency and/or decrease in risk) comparing the portfolios to single design base cases. Hence we first obtain efficiency $\varepsilon^*$ and risk $\mathcal{R}^*$ of a base design, which is calculated though 120 simulations at 5%, 10% and 20% uncertainty levels. The base case is a nominally typical design with 0.2 m$^3$/s injection rate, 50 m stage length, 0.003 Pa.s fluid viscosity, and uniform spacing of 0.5. The perforation factor (Equation S14 detailed in Supplementary Section S3.1) is set as $1.06*10^{10}$ Pa.s$^2$/m$^6$ to ensure the pressure of entry loss is around $1.38*10^7$ Pa (~2000 psi).

Upon establishing a base case, an optimal, non-uniform spacing (NFS) high efficiency and an optimal low risk (i.e. EXL) portfolio are obtained for each case by the previously described methods. Here 6 designs are used for each portfolio, with the details of these designs given in Tables 1 and 2, and the proportion of each portfolio in the optimal designs given in Figure 3(b) (both in tabular form and represented in a spider plot). For each of these portfolios, the efficiency and relative risk is computed and presented as a value normalized by the efficiency and relative risk for the base, single design case. Hence, for each uncertainty level, there is a high efficiency portfolio located at right top and low risk portfolio locates left bottom in Figure 3(b). The tangent line from right to left is again representing a mixture of high efficiency and low risk portfolios. The location in this normalized space is therefore representative of the increase in efficiency (ratio greater than 1) and decrease in relative risk (ratio less than 1) obtained by portfolio design compared to the base case.

Figure 3(b) firstly shows that the efficiency ratio is greater than 1 and relative risk ratio is less than 1 for all cases, indicating that portfolios improve predicted outcomes compared to single repeated designs. In more detail, the efficiency improvement for low risk (EXL) designs are less impacted than high efficiency designs, which is to be expected because the EXL designs ensure that efficiency will remain low due to the high perforation friction. However, the risk reduction for portfolio EXL designs compared to a single base case is profound – around a factor of 5 (relative risk ratio around 0.2) - for all cases.

On the other hand, Figure 3(b) indicates a substantial, 3-4 fold improvement in efficiency for the high efficiency (low perforation friction) designs when a portfolio is used rather than a single design. The risk is also reduced for these cases, with the greatest risk reduction (factor of ~4) corresponding, as expected, to the 20% variability example and a still substantial but lowest (factor of ~2) reduction for the 5% variability example. However, because of this relatively low reduction in risk for the 5% variability cases, the benefit of mixing with EXL cases is much higher, as evidenced by the lower slope (less loss of efficiency for each incremental reduction in risk) of the line connecting high efficiency to EXL portfolios for the 5% case compared to the 10% and 20% cases. So, in summary, the results show potential for 3-5 fold improvements in efficiency and risk can be attained by using portfolios rather than single designs, with the largest benefits of diversification in general coming for the reservoirs with largest variability, and the largest benefits of mixing high efficiency with EXL portfolios coming for the lower risk cases.

**Table 1. Non-Uniform Portfolio Designs at 5%, 10% and 20% Uncertainty**

| Uncertain Level | 5% | | | | | |
|---|---|---|---|---|---|---|
| Design | 1 | 2 | 3 | 4 | 5 | 6 |
| $h_1/Z$ (spacing ratio) | 0.36 | 0.35 | 0.35 | 0.36 | 0.35 | 0.36 |
| $\mu$ (viscosity Pa.s) | 0.54 | 0.21 | 0.17 | 0.21 | 0.20 | 0.57 |
| $T$ (treating time s) | $2.9*10^3$ | $3.6*10^3$ | $2.4*10^3$ | $3.8*10^3$ | $2.6*10^3$ | $2.6*10^3$ |
| $Z$ (stage length m) | 21 | 30 | 28 | 31 | 29 | 20 |
| $Q_o$ (injection rate m$^3$/s) | 0.10 | 0.13 | 0.17 | 0.12 | 0.16 | 0.11 |
| $p_{perf}$ (perforation factor Pa.s$^2$/m$^6$) | $1.3*10^9$ | $3.1*10^6$ | $1.5*10^6$ | $1.2*10^7$ | $6.6*10^5$ | $8.4*10^8$ |
| Uncertain Level | 10% | | | | | |
| Design | 1 | 2 | 3 | 4 | 5 | 6 |
| $h_1/Z$ (spacing ratio) | 0.35 | 0.34 | 0.35 | 0.36 | 0.35 | 0.35 |
| $\mu$ (viscosity Pa.s) | 0.18 | 0.43 | 0.17 | 0.21 | 0.20 | 0.21 |
| $T$ (treating time s) | $2.3*10^3$ | $4.1*10^3$ | $2.4*10^3$ | $3.8*10^3$ | $2.6*10^3$ | $3.6*10^3$ |
| $Z$ (stage length m) | 25 | 30 | 28 | 31 | 29 | 30 |
| $Q_o$ (injection rate m$^3$/s) | 0.15 | 0.11 | 0.17 | 0.12 | 0.16 | 0.13 |
| $p_{perf}$ (perforation factor Pa.s$^2$/m$^6$) | $1.4*10^6$ | $1.7*10^7$ | $1.5*10^6$ | $1.2*10^7$ | $6.6*10^5$ | $3.1*10^6$ |
| Uncertain Level | 20% | | | | | |
| Design | 1 | 2 | 3 | 4 | 5 | 6 |
| $h_1/Z$ (spacing ratio) | 0.29 | 0.35 | 0.37 | 0.36 | 0.34 | 0.36 |
| $\mu$ (viscosity Pa.s) | 0.50 | 0.17 | 0.22 | 0.21 | 0.32 | 0.23 |
| $T$ (treating time s) | $1.5*10^3$ | $2.4*10^3$ | $2.1*10^3$ | $3.8*10^3$ | $2.3*10^3$ | $2.1*10^3$ |
| $Z$ (stage length m) | 24 | 28 | 24 | 31 | 27 | 25 |
| $Q_o$ (injection rate m$^3$/s) | 0.24 | 0.17 | 0.17 | 0.12 | 0.18 | 0.18 |
| $p_{perf}$ (perforation factor Pa.s$^2$/m$^6$) | $1.4*10^5$ | $1.5*10^6$ | $1.3*10^5$ | $1.2*10^7$ | $1.1*10^5$ | $6.8*10^5$ |

**Table 2. EXL Portfolio Designs at 5%, 10% and 20% Uncertainty**

| Uncertain Level | 5% | | | | | |
|---|---|---|---|---|---|---|
| Design | A | B | C | D | E | F |
| $p_{perf}$ (perforation factor $Pa.s^2/m^6$) | $1.3*10^9$ | $3.1*10^6$ | $1.5*10^6$ | $1.2*10^7$ | $6.6*10^5$ | $8.4*10^8$ |
| Uncertain Level | 10% | | | | | |
| Design | A | B | C | D | E | F |
| $p_{perf}$ (perforation factor $Pa.s^2/m^6$) | $1.4*10^6$ | $1.7*10^7$ | $1.5*10^6$ | $1.2*10^7$ | $6.6*10^5$ | $3.1*10^6$ |
| Uncertain Level | 20% | | | | | |
| Design | A | B | C | D | E | F |
| $p_{perf}$ (perforation factor $Pa.s^2/m^6$) | $1.4*10^5$ | $1.5*10^6$ | $1.3*10^5$ | $1.2*10^7$ | $1.1*10^5$ | $6.8*10^5$ |

Besides the improvements available from portfolio combinations in general, there are important observations that can be made in the details of the optimal designs themselves (Table 1). Firstly, we see that in all cases, an optimal non-uniform spacing results when the spacing between fractures 1-2 is about 0.5 times the spacing between fracture 2-3 (and by symmetry, 4-5 compared to 3-4, see Figure 2(c)). There is somewhat less variation in the optimal design values of this non-uniform spacing at the 5% uncertainty level compared to 10% and 20%. Meanwhile, stage length varies from 20-31 meters, without apparent systematic impact from the uncertainty level.

The design details in Table 1 also show that the viscosity for the optimal designs varies between 0.17 to 0.57 Pa.s, which is close to linear gel viscosity. It is predicted by simulations (16) that higher viscosity can lead to more uniform distribution of fluid, which potentially explains why lower viscosity cases would have been likely to give higher relative risk and therefore did not occupy the Pareto set for these examples.

Finally, and perhaps most strikingly, the injection rate and perforation pressure are observed to have a substantial impact on risk but very little impact on the expected value of the efficiency for each design, that is, as long as perforation pressure is small enough to be dominated by other factors impacting the pumping power requirements. In fact, because of this relative insensitivity of the efficiency to these two parameters, they are is observed to vary more widely than other parameters from one design to another within each portfolio described by Table 1. If one were to consider fracture designs deterministically, and therefore independent of variability of the reservoir properties, the conclusion would be that these two parameters have little impact. However, it is shown here that portfolio combinations, with differences between individual designs mostly coming from these "insensitive" parameters, provide a diversity that substantially decreases risk. Once again, the benefit is seen to increase with increasing levels of reservoir variability, i.e. from 5% to 20%. However, it is also important to realize that sensitivity to the perforation pressure loss begins to diminish as the variability increases to 20%, which is indicative of a point where variability of rock properties is so strong that, even if the perforation pressure is large enough to ensure nearly uniform distribution of fluid flow to each fracture, the fracture growth would still be substantially different. At some point the variability of fracture growth is inevitable and there is little systematic benefit from increasing perforation friction, even in EXL designs, a result that is in contrast to the trend that higher perforation pressure should always be used in strongly heterogeneous formations.

## 4. Discussion and Conclusion

The variability and uncertainty of reservoir properties leads to variable and uncertain outcomes of hydraulic fracturing that is essential for oil and gas well stimulation. Such variability leads to both risk and inefficiency, with the latter tied to general trends of stimulation designs that require ever-growing treating pressures, quantities of materials, and durations of pumping. These general trends drive an increase not only in economic inefficiency, but also in environmental impacts associated with extraction of each unit of oil or gas that includes emissions associated with transport of large quantities of materials, risks of inadvertent releases to surface or shallow groundwater, and injection-induced seismicity (22-24). In contrast, we here propose a new concept and method in risk reduction that reverses these troubling trends because it does not require decreasing efficiency, and in fact in many cases can improve efficiency while reducing risk.

The novel concept is developing diverse portfolios comprised of several designs, mixed so as to optimize the risk-return relationship across multiple wells in a given reservoir development.

The approach is applicable not only to hydraulic fracturing, but to a class of design problems with continuously-definable efficiency (return) and risk that is associated with spatio-temporally heterogeneous uncertainty. For the specific case of hydraulic fracturing, by adopting methods for seeking optimal portfolios of investments to rely on metrics definable within mechanical models of hydraulic fracture growth, we show that a portfolio of designs can provide increases in efficiency up to 3.8 times and decreases in risk up to 83% compared to single designs. We also show that the benefits of design portfolios increase with increasing variability of reservoir properties, that is, the more variability and uncertainty one encounters in rock properties and in situ stresses, the more it will be expected that a portfolio of designs will provide higher efficiency (i.e. production per unit of energy input) and lower risk (i.e. more uniformity of stimulation along each well and more similarity of production from well-to-well).

Is important to recognize that the pursuit of a portfolio of designs is not the same as pursuing tailored designs to a variety of known reservoir situations. While the latter approach relies on reducing uncertainty through (potentially costly) additional reservoir characterization and then deterministic application of simulators (and/or experience) to drive customized designs, our approach recognizes that uncertainty is inevitable and therefore should be accounted for in the optimization method. When uncertainty is considered, the results show that design parameters that would have been relatively unimportant to optimizing designs a deterministic framework (i.e. low model sensitivity, in these examples to pumping rate and perforation friction) become important for reduction of variability. Hence the results show that uncertainty fundamentally alters the nature of the optimization problem for oil and gas wells.

Going forward, uncertainty and variability of rock properties is ubiquitous and inevitable in oil and gas development and production. However, from the results presented here it is promising that such variability need not be addressed through methods that drive ever-decreasing efficiency by entailing a drive to ever-increasing injection pressures and/or volumes. Instead, the focus should shift away from the pursuit of single, "optimal" designs to be applied in a factory-style development and production mode. Instead, the industry should pursue model (and perhaps also experience) based identification of multiple designs, each of which is efficiently effective on its own, and then mixing those designs in optimal proportions to attain reduction in risk that is unachievable by any single design that is repetitively and exclusively applied. Furthermore, the potential benefit of diverse design portfolios is applicable not just to oil/gas well stimulation, but to a class of design problems with continuous efficiency and heterogeneous uncertainty, including renewable energy, agriculture, and epidemiology. Hence, this work demonstrates a promising path wherein pursuit of design diversity can combat the negative impacts of variability and uncertainty in oil and gas well stimulation and beyond.

The approach of developing diverse portfolios of designs is shown here to have significant potential to reduce outcome uncertainty for oil and gas well stimulation. However, the essence of the problem points to potential for application in a wide range of fields. Here we firstly discuss the essential characteristics of oil and gas well stimulation that make it suited to a design portfolio approach, after which will discuss other areas of science and engineering that encounter design problems of an ostensibly similar nature.

The first essential characteristic of the problem is that oil and gas stimulation encounters uncertain heterogeneity of the rock formation. Hence, outcomes are impacted not only by the rock property variation range (which is a globally-definable variable, i.e. a mean and standard deviation of each property), but also the fact that its uncertainty is heterogeneous. In this case the heterogeneity is of a spatial nature, meaning that each well within a region draws a set of uncertain properties that, while following overall distributions (the global variables), are essentially independent of the properties of other wells. Of course, there is some spatial correlation to rock properties, however, such correlations are usually at either basin scale (i.e. tendencies that systematically vary across tens to hundreds of kilometers) or rock fabric scale (i.e. tendencies that correlate reservoir rock properties within centimeters to tens of meters). At the well-to-well scale within a region (kilometers to tens of kilometers), the uncertainty associated with heterogeneity is essentially random. Hence, the performance of a given design on one well is not correlated to the performance of another design on another well, meaning that the full benefit of diverse designs can be realized (recalling the reduction in benefit that would arrive from covariance, see Equation 9). Portfolio design diversity in the manifestation presented in this paper is applicable only if the system possesses uncertain heterogeneity.

A second, and perhaps even more important characteristic of oil and gas stimulation that makes portfolio design diversity advantageous is the fact that desirability (in this case efficiency) of outcomes is continuous. This stands in contrast to cases where the outcome is (at least more less) binary – i.e. survival or failure of a structure.

A variety of design problems satisfy the criterial of uncertain heterogeneity and continuously definable desirability of outcomes. Examples include:

1) In Renewable Energy Engineering, the layout of a wind farm, where diversity in design of of blade size/shape, tower spacing, etc. could potentially reduce uncertainty of outcomes such as the (continuously-definable) energy return on energy invested for a system. Here the uncertainty of outcomes arises from both spatially and temporally heterogeneous uncertainty of wind velocity (28).

2) In Epidemiology, where diversity in design for the roll out of a vaccine with continuously-definable efficiency can be defined as decrease in infection rate per vaccine administered. Design parameters include dosing amounts and frequencies, choice of vaccination from multiple available, spacing between priority vaccination areas and timing of administrating the vaccine therein, target proportions of the population to be vaccinated in each area in each phase of the vaccination effort, and so forth. The uncertainty in the outcome is related to heterogeneous variability of each vaccinated individual and of transmission intensities across sectors of the population (29).

Recent examples exist in forestry and agriculture also prove the value of diversity (30,31), although to this point the spatio-temporal heterogeneity of the uncertainty and its implications for optimizing the efficiency-risk relationship has yet to be recognized. Still, from these examples there is a generality that emerges wherein outcome uncertainty is associated with system uncertain heterogeneity and continuously definable. Hence, the framework inspired by MPT and CML, reveals one possible way that uncertain heterogeneity can be managed in a wide range of applications in science and engineering.

## Acknowledgments

This material is based upon work supported by the University of Pittsburgh Center for Energy, Swanson School of Engineering, Department of Chemical and Petroleum Engineering, and Department of Civil and Environmental Engineering. Additional support for recent advances to this work was provided by the National Science Foundation under Grant No. 1645246.